\documentclass[10pt,aps,prb,reqno]{amsart}
\usepackage{amsmath}
\usepackage{amssymb}

\newtheorem{theorem}{Theorem}[section]
\newtheorem{remark}{Remark}[section]
\newtheorem{lemma}[theorem]{Lemma}
\newtheorem{proposition}[theorem]{Proposition}

\begin{document}
\title[Logarithmically supercritical MHD]{Global regularity of logarithmically supercritical MHD system with zero diffusivity}
\author{Kazuo Yamazaki}  
\date{}
\thanks{The author expresses gratitude to Professor Jiahong Wu and Professor David Ullrich for their teaching and reviewers for their comments that improved this manuscript greatly. }
\maketitle 

\begin{abstract}
We prove the global regularity of the solution pair to the N-dimensional logarithmically supercritical magnetohydrodynamics system with zero diffusivity. This is the endpoint case omitted in the work of [24]; it also improves some previous results logarithmically. 

\vspace{5mm}

\textbf{Keywords: Global regularity, magnetohydrodynamics system, Navier-Stokes system}
\end{abstract}
\footnote{2000MSC : 35B65, 35Q35, 35Q86}
\footnote{Department of Mathematics, Oklahoma State University, 401 Mathematical Sciences, Stillwater, OK 74078, USA}

\section{Introduction and statement of results}

We study the N-dimensional generalized magnetohydrodynamics (MHD) system:

\begin{equation}
\begin{cases}
\partial_{t}u + (u\cdot\nabla)u - (b\cdot\nabla)b + \nabla \pi + \nu\mathcal{L}_{1}^{2} u = 0,\\
\partial_{t}b + (u\cdot\nabla)b - (b\cdot\nabla)u + \eta\mathcal{L}_{2}^{2}b = 0,\\
\nabla\cdot u = \nabla\cdot b = 0, \hspace{5mm} (u,b)(x,0) = (u_{0}, b_{0})(x), \hspace{5mm} x \in \mathbb{R}^{N}, N \geq 2,
\end{cases}
\end{equation}
where $u: \mathbb{R}^{N}\times \mathbb{R}^{+}\mapsto \mathbb{R}^{N}$ represents the velocity vector field, $b: \mathbb{R}^{N}\times \mathbb{R}^{+}\mapsto \mathbb{R}^{N}$ the magnetic vector field, $\pi: \mathbb{R}^{N}\times \mathbb{R}^{+}\mapsto \mathbb{R}$ the pressure scalar field and $\nu, \eta \geq 0$ the kinematic viscosity and diffusivity constants respectively. Moreover, we denote operators $\mathcal{L}_{1}, \mathcal{L}_{2}$ defined by 

\begin{equation}
\widehat{\mathcal{L}_{1}f}(\xi) = m_{1}(\xi) \hat{f}(\xi), \hspace{5mm} \widehat{\mathcal{L}_{2}f}(\xi) = m_{2}(\xi)\hat{f}(\xi),
\end{equation}
where $m_{i}(\xi), i = 1, 2$ obeys the lower bound

\begin{equation}
m_{1}(\xi) \geq \frac{\lvert \xi\rvert^{\alpha}}{g_{1}(\xi)}, \hspace{5mm} m_{2}(\xi) \geq \frac{\lvert \xi\rvert^{\beta}}{g_{2}(\xi)}, 
\end{equation}
for all sufficiently large $\lvert \xi\rvert$ and $g_{i}: \mathbb{R}^{+}\mapsto \mathbb{R}^{+}, g_{i} \geq 1, i = 1, 2$ are radially symmetric, non-decreasing functions. 

It is well-known that in case $N = 2, 3, \nu, \eta > 0, \alpha = \beta = 1, g_{i} \equiv 1, i = 1, 2$, the MHD system possesses at least one global $L^{2}$ weak solution for any initial data pair $(u_{0}, b_{0}) \in L^{2}$; in case $N =2$, the solution is unique (cf. [19]). In fact, in any dimension $N \geq 2$, the case $\nu, \eta > 0, \alpha \geq \frac{1}{2} + \frac{N}{4}, \beta \geq \frac{1}{2} + \frac{N}{4}, g_{i} \equiv 1, i = 1, 2$ implies that the system is energy-critical and the existence of the unique global strong solution can be shown (cf. [26] and interesting improvement in [4] and [29] for the case $N=2$).   

In [20] for the wave equation and then in [21] for the Navier-Stokes system (NSE) for $N \geq 3$, the author showed global regularity of the solution even in the logarithmically supercritical regime; a different proof in the case of the NSE appeared recently in [14]. More examples of this new type of result followed: [5] in the case of two-dimensional Euler equation, [12] the two-dimensional Boussinesq system, [18] three-dimensional wave equation again, [27] the three-dimensional Leray-alpha type models. In particular in [24], the case $\nu, \eta > 0, \alpha \geq \frac{1}{2} + \frac{N}{4}, \beta > 0$ such that $\alpha + \beta \geq 1 + \frac{N}{2}$, the system (1) with $g_{i}, i = 1,2$ satisfying certain conditions not necessarily that $g_{i} \equiv 1, i = 1, 2$, global regularity result was shown. 

In general even without the logarithmic supercriticality, the endpoint case $\eta = 0$ requires a more subtle proof as non-linear terms such as $(b\cdot\nabla) b$ forces one to rely on the logarithmic-type inequality from [1] due to the lack of diffusivity. Indeed, the requirement that $\eta > 0, \beta > 0$ is crucial in the work of [24]. 

On the other hand, numerical analysis results (e.g. [10], [17]) indicate more dominant role played by the velocity vector field in preserving the regularity of the solution pair. Moreover, besides regularity criteria that depends on both the velocity and magnetic vector fields that we saw in the past (e.g. [2], [9]), starting from the pioneering works of [11] and [32], recently various regularity criteria of the MHD system in terms of only the velocity vector field appeared (e.g. [3], [8], [25], [28], [31]). Thus, it is of interest whether we can extend the results of [24] to the endpoint case of zero diffusivity. We answer this question: 

\begin{theorem}
Let $\nu > 0, \eta =0, \alpha \geq 1 + \frac{N}{2}$ and $g_{1}: \mathbb{R}^{+} \mapsto \mathbb{R}^{+}$ be a radially symmetric, non-decreasing  function such that $g_{1} \geq 1$ and satisfy 

\begin{equation}
\int_{e}^{\infty}\frac{d\tau}{g_{1}^{2}(\tau)\ln (\tau)\tau}= \infty.
\end{equation}
Then for all initial data pair $(u_{0}, b_{0})\in H^{s}, s \geq 3 + N$, there exists a unique global classical solution pair $(u,b)$ to (1) where $\mathcal{L}_{1}$ is defined by (2) and (3). 

\end{theorem}

\begin{remark}
\begin{enumerate}
\item There are various ways to obtain different initial regularity conditions. We chose the statement above for simplicity; its proof follows the argument in [16] using mollifiers.  
\item Theorem 1.1 completes all cases of global regularity of logarithmically supercritical MHD system with equal or more dependence on the dissipation than diffusion. It also improves the result of [23] logarithmically. 
\item After this work was completed, we were informed of the work in [22]. We note that, as the authors in [22] acknowledge, their proof was inspired by the work in [15] while our proof was largely inspired by the work of [30] and partially [7] and [32]. In particular, inspired by the work of [7], we chose not to apply the Littlewood-Paley decomposition on the system itself but rather only on $\nabla u$, in contrast to the proof of [24]. Naturally both proofs from [22] and this manuscript rely on the Brezis-Wainger type inequality. However, the precise way of splitting $\lVert \nabla u\rVert_{L^{\infty}}$ into two parts is different (see Lemma 2.1 of [22] and (6) of this manuscript). 

Moreover, compared to our condition of (4) on $g_{1}$, the authors in [22] have the condition on $g_{1}$ such that there exists an absolute constant $c \geq 0$ satisfying 
\begin{equation*}
g_{1}^{2}(\tau) \leq c \ln (e+\tau).
\end{equation*}
Therefore, (4) of this manuscript in particular allows $g_{1}(\tau)$ to grow as $\sqrt{\ln(e + \ln(e+\tau))}$ while the result in [22] allows $g_{1}(\tau)$ to grow as $\sqrt{\ln(e+\tau)}$. On the other hand, the integral-type condition of (4) in this manuscript allows $g_{1}(\tau)$ to have spikes from time to time, which is a situation not covered in [22].  
\end{enumerate}

\end{remark}

In the Preliminary section, we set up notations and state key lemmas; thereafter, we prove our theorem. 

\section{Preliminary}

Let us use the notation $A\lesssim_{a, b}B$ to imply that there exists a positive constant $c$ that depends on $a, b$ such that $A \leq cB$. We use the following well-known commutator estimate: 

\begin{lemma}
(cf. [13]) Let $f,g$ be smooth such that $\nabla f \in L^{p_{1}}, \Lambda^{s-1}g \in L^{p_{2}}, \Lambda^{s}f \in L^{p_{3}}, g \in L^{p_{4}}, p \in (1,\infty), p_{2}, p_{3} \in (1,\infty), \frac{1}{p} = \frac{1}{p_{1}}+\frac{1}{p_{2}} = \frac{1}{p_{3}} + \frac{1}{p_{4}}, s > 0.$ Then there exists a constant $c > 0$ such that

\begin{equation*}
\lVert \Lambda^{s}(fg) - f\Lambda^{s}g\rVert_{L^{p}} \leq c(\lVert \nabla f\rVert_{L^{p_{1}}}\lVert \Lambda^{s-1}g\rVert_{L^{p_{2}}} + \lVert \Lambda^{s}f\rVert_{L^{p_{3}}}\lVert g\rVert_{L^{p_{4}}}).
\end{equation*}

\end{lemma} 

Let us recall the notion of Besov spaces (cf. [6]). We denote by $\mathcal{S}(\mathbb{R}^{N})$ the Schwartz class functions and $\mathcal{S}'(\mathbb{R}^{N})$, its dual. We define $\mathcal{S}_{0}$ to be the subspace of $\mathcal{S}$ in the following sense:

\begin{equation*}
\mathcal{S}_{0} = \{\phi \in \mathcal{S}, \int_{\mathbb{R}^{N}}\phi(x)x^{\gamma}dx = 0, \lvert\gamma\rvert = 0, 1, 2, ... \}.
\end{equation*}
Its dual $\mathcal{S}_{0}'$ is given by $\mathcal{S}_{0}' = \mathcal{S}/\mathcal{S}_{0}^{\perp} = \mathcal{S}'/\mathcal{P}$ where $\mathcal{P}$ is the space of polynomials. For $j \in \mathbb{Z}$ we define 

\begin{equation*}
A_{j} = \{\xi\in \mathbb{R}^{N}: 2^{j-1} < \lvert\xi\rvert < 2^{j+1}\}.
\end{equation*}
It is well-known that there exists a sequence $\{\Phi_{j}\} \in \mathcal{S}(\mathbb{R}^{N})$ such that 

\begin{equation*}
\text{ supp }\hat{\Phi}_{j}\subset A_{j}, \hspace{5mm} \hat{\Phi}_{j}(\xi) = \hat{\Phi}_{0}(2^{-j}\xi) \hspace{5mm} \text{or} \hspace{5mm} \Phi_{j}(x) = 2^{jN}\Phi_{0}(2^{j}x) \hspace{5mm} \text{ and } 
\end{equation*}

\begin{equation*}
\sum_{k = -\infty}^{\infty}\hat{\Phi}_{k}(\xi) = 
\begin{cases}
1 \hspace{2mm} \text{ if } \hspace{1mm} \xi \in \mathbb{R}^{N}\setminus \{0\},\\
0 \hspace{2mm} \text{ if } \hspace{1mm} \xi = 0.
\end{cases}
\end{equation*}
Consequently, for any $f \in \mathcal{S}_{0}'$, 

\begin{equation*}
\sum_{k = -\infty}^{\infty}\Phi_{k}\ast f = f.
\end{equation*}

To define the homogeneous Besov space, we set 

\begin{equation*}
\triangle_{j}f = \Phi_{j}\ast f, \hspace{5mm} j = 0, \pm 1, \pm 2, ....
\end{equation*}
With such we can define for $s \in \mathbb{R}, p, q \in [1,\infty]$, the homogeneous Besov space

\begin{equation*}
\dot{B}_{p, q}^{s} = \{f \in \mathcal{S}_{0}' : \lVert f\rVert_{\dot{B}_{p, q}^{s}} < \infty\},
\end{equation*}
where 

\begin{equation*}
\lVert f\rVert_{\dot{B}_{p,q}^{s}} = 
\begin{cases}
(\sum_{j}(2^{js}\lVert\triangle_{j}f\rVert_{L^{p}})^{q})^{\frac{1}{q}} & \text{if } q < \infty,\\
\sup_{j}2^{js}\lVert\triangle_{j}f\rVert_{L^{p}} & \text{if } q = \infty.
\end{cases}
\end{equation*}

To define the inhomogeneous Besov space, we let $\Psi \in C_{0}^{\infty}(\mathbb{R}^{N})$ be such that

\begin{equation*}
1 = \hat{\Psi}(\xi) + \sum_{k=0}^{\infty}\hat{\Phi}_{k}(\xi), \hspace{5mm} \Psi\ast f + \sum_{k=0}^{\infty}\Phi_{k}\ast f = f,
\end{equation*}
for any $f \in \mathcal{S}'$. With that, we set

\begin{equation*}
\triangle_{j}f = 
\begin{cases}
0 &\text{ if } j \leq -2,\\
\Psi\ast f &\text{ if } j = -1,\\
\Phi_{j}\ast f &\text{ if } j = 0, 1, 2, ...,
\end{cases}
\end{equation*}
and define for any $s \in \mathbb{R}, p, q \in [1, \infty]$, the  inhomogeneous Besov space 

\begin{equation*}
B_{p, q}^{s} = \{f \in \mathcal{S}' : \lVert f\rVert_{B_{p, q}^{s}} < \infty \},
\end{equation*}
where

\begin{equation*}
\lVert f\rVert_{B_{p, q}^{s}} = 
\begin{cases}
(\sum_{j = -1}^{\infty}(2^{js}\lVert\triangle_{j}f\rVert_{L^{p}})^{q})^{\frac{1}{q}}, &\text{if } q < \infty,\\
\sup_{-1 \leq j < \infty}2^{js}\lVert \triangle_{j}f\rVert_{L^{p}} &\text{if } q = \infty.
\end{cases}
\end{equation*}
In particular $B_{2,2}^{s} = H^{s}$. 

Finally, the following lemma will be useful:

\begin{lemma} (cf. [6])
Bernstein's Inequality: Let f $\in L^{p}(\mathbb{R}^{N})$ with 1 $\leq p \leq q \leq \infty$ and $0 < r < R$. Then for all $k \in \mathbb{Z}^{+}\cup\{0\}$, and $\lambda > 0$, there exists a constant $C_{k} > 0$ such that

\begin{equation*}
\begin{cases}
\sup_{\lvert \gamma \rvert = k} \lVert\partial^{\gamma}f\rVert_{L^{q}} \leq C_{k}\lambda^{k+N(\frac{1}{p} - \frac{1}{q})}\lVert f\rVert_{L^{p}}  & \text{ if } \text{supp }\hat{f} \subset \{\xi : \lvert\xi\rvert \leq \lambda r\},\\
C_{k}^{-1}\lambda^{k}\lVert f \rVert_{L^{p}} \leq sup_{\lvert\gamma\rvert = k} \lVert\partial^{\gamma}f\rVert_{L^{p}} \leq C_{k}\lambda^{k}\lVert f \rVert_{L^{p}}  & \text{ if  supp } \hat{f} \subset \{\xi : \lambda r \leq \lvert\xi\rvert \leq \lambda R\},
\end{cases}
\end{equation*}
and if we replace derivative $\partial^{\gamma}$ by the fractional derivative, the inequalities remain valid only with trivial modifications.
\end{lemma}

\section{Proof}

Let $s$ be a large integer. To establish global regularity, it suffices to prove an \textit{a priori } bound of the form 

\begin{equation*}
\lVert u(t)\rVert_{H^{s}} \leq C(s, \lVert u_{0}\rVert_{H^{s}}, \lVert b_{0}\rVert_{H^{s}}, T, g_{1}),
\end{equation*}
for all $0 \leq t \leq T < \infty$ where $C(s, \lVert u_{0}\rVert_{H^{s}}, \lVert b_{0}\rVert_{H^{s}}, T, g_{1})$ is a constant depending on $s, \lVert u_{0}\rVert_{H^{s}}, \lVert b_{0}\rVert_{H^{s}}, T$ and $g_{1}$. We now fix $(u_{0}, b_{0}), (u, b), T$ and assume without loss of generality that $\nu = 1$. Taking inner products of the first and second equation with $u$ and $b$ respectively, we obtain due to the incompressibility of $u$ and $b$

\begin{eqnarray*}
\frac{1}{2}\partial_{t}(\lVert u\rVert_{L^{2}}^{2} + \lVert b\rVert_{L^{2}}^{2}) + \lVert \mathcal{L}_{1}u\rVert_{L^{2}}^{2} = 0.
\end{eqnarray*}

Therefore,

\begin{equation}
\sup_{t\in [0,T]}\lVert u(t)\rVert_{L^{2}}^{2} + \lVert b(t)\rVert_{L^{2}}^{2} + 2\int_{0}^{T}\lVert \mathcal{L}_{1}u\rVert_{L^{2}}^{2} d\tau \lesssim_{u_{0}, b_{0}, T} 1.
\end{equation}

Now let us denote by 

\begin{equation*}
X(t) := \lVert \nabla u(t)\rVert_{L^{2}}^{2} + \lVert \nabla b(t)\rVert_{L^{2}}^{2}, \hspace{5mm} Y(t) := \lVert \Lambda^{s}u(t)\rVert_{L^{2}}^{2} + \lVert \Lambda^{s}b(t)\rVert_{L^{2}}^{2}.
\end{equation*}

\subsection{$H^{1}$-estimate}

We prove the following proposition:

\begin{proposition}

Let $\nu > 0, \eta =0, \alpha \geq 1 + \frac{N}{2}$ and $g_{1}: \mathbb{R}^{+} \mapsto \mathbb{R}^{+}$ be a radially symmetric, non-decreasing function such that $g_{1} \geq 1$ and satisfy (4). Then the solution pair $(u,b)$ to (1) in $[0,T]$ satisfies 

\begin{equation*}
\sup_{t\in [0,T]} \lVert \nabla u(t)\rVert_{L^{2}}^{2} + \lVert \nabla b(t)\rVert_{L^{2}}^{2} + \int_{0}^{T}\lVert \mathcal{L}_{1}\nabla u\rVert_{L^{2}}^{2} d\tau < \infty.
\end{equation*}

\end{proposition}

We apply $\nabla$ on the first and second equations of (1), take $L^{2}$-inner products with $\nabla u$ and $\nabla b$ respectively and sum to obtain 

\begin{eqnarray*}
\frac{1}{2}\partial_{t}X(t) + \lVert \mathcal{L}_{1}\nabla u\rVert_{L^{2}}^{2}  \lesssim \int \lvert \nabla u\rvert^{3} + \lvert \nabla u\rvert \lvert\nabla b\rvert^{2}.
\end{eqnarray*}

We apply H$\ddot{o}$lder's inequality to obtain

\begin{eqnarray*}
\frac{1}{2}\partial_{t}X(t) + \lVert \mathcal{L}_{1}\nabla u\rVert_{L^{2}}^{2} \lesssim \lVert \nabla u\rVert_{L^{\infty}}(\lVert \nabla u\rVert_{L^{2}} + \lVert \nabla b\rVert_{L^{2}})^{2}.
\end{eqnarray*}

By Bernstein's inequality for some $M_{1} > 0$ to be determined subsequently we obtain

\begin{eqnarray*}
\lVert \nabla u\rVert_{L^{\infty}} &\leq& \sum_{j\geq -1}\lVert \Delta_{j}\nabla u\rVert_{L^{\infty}}\\
&\lesssim&  \sum_{2^{j} \leq M_{1}}\frac{2^{j(\frac{N}{2})}}{g_{1}(2^{j})}\lVert \Delta_{j}\nabla u\rVert_{L^{2}}g_{1}(2^{j}) + \sum_{2^{j} > M_{1}}2^{-j}\frac{2^{j(\frac{N}{2} + 1)}}{g_{1}(2^{j})}\lVert \Delta_{j}\nabla u\rVert_{L^{2}}g_{1}(2^{j}).
\end{eqnarray*}

Since $g_{1}$ is increasing, by (3) we obtain

\begin{eqnarray*}
\lVert \nabla u\rVert_{L^{\infty}}
\lesssim g_{1}(M_{1})\sum_{2^{j} \leq M_{1}}\lVert \Delta_{j}\mathcal{L}_{1}u\rVert_{L^{2}} + \sum_{2^{j} > M_{1}}2^{-j}g_{1}(2^{j})\lVert \Delta_{j}\mathcal{L}_{1}\nabla u\rVert_{L^{2}}.
\end{eqnarray*}

We further bound using H$\ddot{o}$lder's inequalities to obtain

\begin{equation}
\lVert \nabla u\rVert_{L^{\infty}} \lesssim g_{1}(M_{1})\sqrt{\ln(M_{1})}\lVert \mathcal{L}_{1}u\rVert_{L^{2}} + M_{1}^{-\frac{1}{2}} \lVert\mathcal{L}_{1}\nabla u\rVert_{L^{2}}.
\end{equation}

Now set $M_{1}:= e + X(t)$ so that by (6) and Young's inequality we have 

\begin{eqnarray*}
&&\lVert \nabla u\rVert_{L^{\infty}} (\lVert \nabla u\rVert_{L^{2}} + \lVert \nabla b\rVert_{L^{2}})^{2}\\
&\lesssim& g_{1}(e+X(t))\sqrt{\ln(e+X(t))}\lVert \mathcal{L}_{1}u\rVert_{L^{2}} X(t) + \lVert \mathcal{L}_{1}\nabla u\rVert_{L^{2}}(\lVert \nabla u\rVert_{L^{2}} + \lVert \nabla b\rVert_{L^{2}})\\
&\leq& \frac{1}{2}\lVert \mathcal{L}_{1}\nabla u\rVert_{L^{2}}^{2} + c\left(g_{1}(e + X(t)\right) \sqrt{\ln(e + X(t))}\lVert \mathcal{L}_{1}u\rVert_{L^{2}}X(t) + cX(t).
\end{eqnarray*}

Absorbing the dissipative term, we have by Young's inequalities again

\begin{eqnarray}
\partial_{t}X(t) + \lVert \mathcal{L}_{1}\nabla u\rVert_{L^{2}}^{2} \lesssim \left(g_{1}^{2}(e + X(t))\ln(e + X(t))\right)(e + X(t))\left(1+\lVert \mathcal{L}_{1}u\rVert_{L^{2}}^{2}\right).
\end{eqnarray}

Thus, for any $t \in [0,T]$, we have due to (5)

\begin{eqnarray*}
\int_{e + X(0)}^{e + X(t)}\frac{d\tau}{g_{1}^{2}(\tau)\ln(\tau)\tau} \lesssim \int_{0}^{T}1+\lVert \mathcal{L}_{1}u\rVert_{L^{2}}^{2} d\tau \lesssim_{u_{0}, b_{0}, T}1.
\end{eqnarray*}

By hypothesis (4) this implies 

\begin{equation}
\sup_{t\in [0,T]}X(t) < \infty.
\end{equation}

From (7) it now follows using (8) that 

\begin{eqnarray}
\int_{0}^{T}\lVert \mathcal{L}_{1}\nabla u\rVert_{L^{2}}^{2} d\tau \lesssim_{u_{0}, b_{0}, g_{1}, T}1.
\end{eqnarray}

This completes the proof of Proposition 3.1. 

\textit{Proof of Theorem 1.1}

We fix 

\begin{equation}
\gamma \in \left(1+\frac{N}{2}, 2+\frac{N}{2}\right),
\end{equation}
apply $\Lambda^{\gamma}$ on both equations of (1), take $L^{2}$-inner products with $\Lambda^{\gamma}u$ and $\Lambda^{\gamma}b$ on the first and second equations respectively and estimate in sum

\begin{eqnarray*}
&&\partial_{t}(\lVert \Lambda^{\gamma}u\rVert_{L^{2}}^{2} + \lVert \Lambda^{\gamma}b\rVert_{L^{2}}^{2}) + \lVert \mathcal{L}_{1}\Lambda^{\gamma}u\rVert_{L^{2}}^{2}\\
&=& -\int [\Lambda^{\gamma}((u\cdot\nabla)u) - u\cdot\nabla\Lambda^{\gamma}u]\cdot\Lambda^{\gamma}u -\int [\Lambda^{\gamma}((u\cdot\nabla)b) - u\cdot\nabla\Lambda^{\gamma}b]\cdot\Lambda^{\gamma}b\\
&&+\int[\Lambda^{\gamma}((b\cdot\nabla)b) - b\cdot\nabla\Lambda^{\gamma}b]\cdot\Lambda^{\gamma}u +\int[\Lambda^{\gamma}((b\cdot\nabla)u) - b\cdot\nabla\Lambda^{\gamma}u]\cdot\Lambda^{\gamma}b.
\end{eqnarray*} 

We bound as follows:

\begin{eqnarray}
&&\partial_{t}(\lVert \Lambda^{\gamma}u\rVert_{L^{2}}^{2} + \lVert \Lambda^{\gamma}b\rVert_{L^{2}}^{2}) + \lVert \mathcal{L}_{1}\Lambda^{\gamma}u\rVert_{L^{2}}^{2}\\
&\lesssim& (\lVert \nabla u\rVert_{L^{\infty}}\lVert \Lambda^{\gamma-1}\nabla u\rVert_{L^{2}} + \lVert\Lambda^{\gamma}u\rVert_{L^{2}}\lVert \nabla u\rVert_{L^{\infty}})\lVert \Lambda^{\gamma}u\rVert_{L^{2}}\nonumber\\
&&+(\lVert \nabla u\rVert_{L^{\infty}} \lVert \Lambda^{\gamma-1}\nabla b\rVert_{L^{2}} + \lVert \Lambda^{\gamma}u\rVert_{L^{2}}\lVert \nabla b\rVert_{L^{\infty}}) \lVert \Lambda^{\gamma}b\rVert_{L^{2}}\nonumber\\
&&+(\lVert \nabla b\rVert_{L^{\infty}}\lVert\Lambda^{\gamma-1}\nabla b\rVert_{L^{2}} + \lVert \Lambda^{\gamma}b\rVert_{L^{2}}\lVert \nabla b\rVert_{L^{\infty}}) \lVert \Lambda^{\gamma}u\rVert_{L^{2}}\nonumber\\
&&+(\lVert \nabla b\rVert_{L^{\infty}}\lVert \Lambda^{\gamma-1}\nabla u\rVert_{L^{2}} + \lVert \Lambda^{\gamma}b\rVert_{L^{2}} \lVert \nabla u\rVert_{L^{\infty}})\lVert \Lambda^{\gamma}b\rVert_{L^{2}}\nonumber\\
&\lesssim& \lVert \nabla u\rVert_{L^{\infty}}(\lVert \Lambda^{\gamma} u\rVert_{L^{2}}^{2} + \lVert \Lambda^{\gamma} b\rVert_{L^{2}}^{2}) + \lVert \Lambda^{\gamma}u\rVert_{L^{2}}\lVert \nabla b\rVert_{L^{\infty}}\lVert\Lambda^{\gamma} b\rVert_{L^{2}},\nonumber
\end{eqnarray}
by H$\ddot{o}$lder's inequalities and Lemma 2.1. By (6) and (8), we obtain an estimate of 

\begin{eqnarray*}
\lVert \nabla u\rVert_{L^{\infty}} &&\lesssim g_{1}(M_{1})\sqrt{\ln(M_{1})}\lVert \mathcal{L}_{1}u\rVert_{L^{2}} + M_{1}^{-\frac{1}{2}} \lVert\mathcal{L}_{1}\nabla u\rVert_{L^{2}}\\
&&\lesssim_{g_{1}} \lVert \mathcal{L}_{1}u\rVert_{L^{2}} + \lVert \mathcal{L}_{1}\nabla u\rVert_{L^{2}}.
\end{eqnarray*}

Using this and a Gagliarido-Nirenberg inequality justified by (10), we continue the estimate of (11) by 

\begin{eqnarray*}
&&\partial_{t}(\lVert \Lambda^{\gamma}u\rVert_{L^{2}}^{2} + \lVert \Lambda^{\gamma}b\rVert_{L^{2}}^{2}) + \lVert \mathcal{L}_{1}\Lambda^{\gamma}u\rVert_{L^{2}}^{2}\\
&\lesssim& (\lVert \mathcal{L}_{1} u\rVert_{L^{2}} + \lVert \mathcal{L}_{1} \nabla u\rVert_{L^{2}}  ) ( \lVert \Lambda^{\gamma} u\rVert_{L^{2}}^{2} + \lVert \Lambda^{\gamma} b\rVert_{L^{2}}^{2})\\
&&+ \lVert \Lambda^{\gamma} u\rVert_{L^{2}}\lVert\nabla b\rVert_{L^{2}}^{\frac{2(\gamma-1)-N}{2(\gamma-1)}}\lVert\Lambda^{\gamma-1} \nabla b\rVert_{L^{2}}^{\frac{N}{2(\gamma-1)}}\lVert\Lambda^{\gamma} b\rVert_{L^{2}}.
\end{eqnarray*}

On the other hand, by a similar procedure as before,

\begin{eqnarray*}
\lVert \Lambda^{\gamma} u\rVert_{L^{2}}
&\leq& \sum_{j\geq -1} \lVert \Lambda^{\gamma} \Delta_{j} u\rVert_{L^{2}}\\
&\lesssim& \sum_{j\geq -1} \lVert \mathcal{L}_{1}\nabla u\rVert_{L^{2}} g_{1}(2^{j}) 2^{j(\gamma-(2+\frac{N}{2}))}
\lesssim \lVert \mathcal{L}_{1}\nabla u\rVert_{L^{2}},
\end{eqnarray*}
due to Bernstein's inequality, (3), H$\ddot{o}$lder's inequality and the fact that $g_{1}$ grows logarithmically while $\gamma < 2+\frac{N}{2}$. Therefore, by Proposition 3.1 and Young's inequality justified by (10) we obtain

\begin{eqnarray*}
&&\partial_{t}(\lVert \Lambda^{\gamma} u\rVert_{L^{2}}^{2} + \lVert \Lambda^{\gamma} b\rVert_{L^{2}}^{2}) + \lVert \mathcal{L}_{1}\Lambda^{\gamma} u\rVert_{L^{2}}^{2}\\
&\lesssim& (\lVert \mathcal{L}_{1} u\rVert_{L^{2}} + \lVert \mathcal{L}_{1}\nabla u\rVert_{L^{2}})(\lVert \Lambda^{\gamma} u\rVert_{L^{2}}^{2} + \lVert \Lambda^{\gamma} b\rVert_{L^{2}}^{2})\\
&&+ \lVert \mathcal{L}_{1}\nabla u\rVert_{L^{2}} \lVert \Lambda^{\gamma} b\rVert_{L^{2}}^{1+\frac{N}{2(\gamma-1)}}\\
&\lesssim& (\lVert \mathcal{L}_{1} u\rVert_{L^{2}} + \lVert \mathcal{L}_{1}\nabla u\rVert_{L^{2}})(\lVert \Lambda^{\gamma} u\rVert_{L^{2}}^{2} + \lVert \Lambda^{\gamma} b\rVert_{L^{2}}^{2} + e).
\end{eqnarray*}
That is, 

\begin{equation*}
\partial_{t} \ln (e+ \lVert \Lambda^{\gamma} u\rVert_{L^{2}}^{2} + \lVert \Lambda^{\gamma} b\rVert_{L^{2}}^{2}) \lesssim \lVert \mathcal{L}_{1} u\rVert_{L^{2}} + \lVert \mathcal{L}_{1}\nabla u\rVert_{L^{2}}.
\end{equation*}

Thus, by (5) and Proposition 3.1, this leads to 

\begin{equation}
\sup_{t\in [0,T]} \lVert \Lambda^{\gamma} u(t)\rVert_{L^{2}}^{2} + \lVert \Lambda^{\gamma}b(t)\rVert_{L^{2}}^{2} + \int_{0}^{T} \lVert \mathcal{L}_{1}\Lambda^{\gamma}u\rVert_{L^{2}}^{2} d\tau \lesssim_{u_{0}, b_{0}, g_{1}, T}1,
\end{equation}
for any $\gamma \in \left(1+ \frac{N}{2}, 2+\frac{N}{2}\right)$. Extending to higher regularity follows from Sobolev inequality: for any $s \in \mathbb{R}^{+}$ this time, by applying Lemma 2.1 as done in (11) we may obtain 

\begin{eqnarray*}
\frac{1}{2}\partial_{t}Y(t) + \lVert\mathcal{L}_{1}\Lambda^{s} u\rVert_{L^{2}}^{2}\lesssim (\lVert \nabla u\rVert_{L^{\infty}} + \lVert\nabla b\rVert_{L^{\infty}})Y(t) \lesssim Y(t),
\end{eqnarray*}
by (10) and (12). This completes the proof of Theorem 1.1.

\end{document}